\documentclass{compositio}
 
\usepackage{bbm, amsmath, amsfonts, amscd, latexsym, amsthm, amssymb, graphicx}
\usepackage[symbol]{footmisc}

\makeatletter
\newif\if@check@engine  \@check@enginetrue 
\makeatother
\usepackage[usenames,dvipsnames]{pstricks}

\newtheorem{theor}{\hspace{1cm}{\sc Theorem}}[section]

\newtheorem{sledst}[theor]{\hspace{1cm}{\sc Corollary}}
\newtheorem{lemma}[theor]{\hspace{1cm}{\sc Lemma}}
\newtheorem{conj}[theor]{\hspace{1cm}{\sc Conjecture}}
\newtheorem*{utver*}{\hspace{1cm}{\sc Proposition}}
\theoremstyle{definition}
\newtheorem{classif}[theor]{\hspace{1cm}{\sc Classification}}
\newtheorem{defin}[theor]{\hspace{1cm}{\sc Definition}}
\newtheorem{exa}[theor]{\hspace{1cm}{\sc Example}}
\newtheorem{rem}[theor]{\hspace{1cm}{\sc Remark}}

\renewcommand{\contentsname}{}

\newcommand{\codim}{\mathop{\rm codim}\nolimits}

\newcommand{\MV}{\mathop{\rm MV}\nolimits}

\newcommand{\supp}{\mathop{\rm supp}\nolimits}
\newcommand{\coker}{\mathop{\rm coker}\nolimits}

\def\R{\mathbb R}

\def\Z{\mathbb Z}

\def\C{\mathbb C}
\def\CC{({\mathbb C}\setminus 0)}

\def\CP{\mathbb C\mathbb P}

\begin{document}

\title[Galois theory for systems of equations]{Galois theory for general systems of polynomial equations}
\author{A. Esterov} 
\email{aesterov@hse.ru}
\address{National Research University Higher School of Economics\newline Faculty of Mathematics NRU HSE, Usacheva str., 6, Moscow, 119048, Russia}
\thanks{Research supported by the Russian Science Foundation grant, project 16-11-10316.}

\begin{abstract} We prove that the monodromy group of a reduced irreducible square system of general polynomial equations equals the symmetric group. This is a natural first step towards the Galois theory of general systems of polynomial equations, because arbitrary systems split into reduced irreducible ones upon monomial changes of variables. 

In particular, our result proves the multivariate version of the Abel--Ruffini theorem: the classification of general systems of equations solvable by radicals reduces to the classification of lattice polytopes of mixed volume 4 (which we prove to be finite in every dimension). We also notice that the monodromy of every general system of equations is either symmetric or imprimitive.

The proof is based on a new result of independent importance regarding dual defectiveness of systems of equations: the discriminant of a reduced irreducible square system of general polynomial equations is a hypersurface unless the system is linear up to a monomial change of variables.
\end{abstract}

\dedication{To R. K. Gordin on the occasion of his 70th birthday.}

\maketitle

\section{Introduction} \label{Sintro}

{\bf Galois theory for lattice polytopes.} A problem of enumerative geometry asks how many geometric objects satisfy a generic geometric constraint in a given space of constraints $P$. Galois theory for this enumerative problem studies how the solutions of this problem permute as the constraint runs along loops in $P$.  In the last decade, particularly strong results were obtained in Galois theory of Schubert calculus, see \cite{schsot} and references therein.

We develop Galois theory in the same vein for another well known enumerative problem: the Kouchnirenko--Bernstein theorem, counting the solutions of a system of generic polynomial equations composed of a given finite collection of monomials. More accurately, let us identify points $a=(a_1,\ldots,a_n)\in\Z^n$ with monomials $x^a=x_1^{a_1}\ldots x_n^{a^n}$, then every finite set of monomials $A\subset\Z^n$ gives rise to the space of Laurent polynomials $\C^A=\{\sum_{a\in A} c_ax^a,\, c_a\in\C\}$, {\it supported at} $A$. These polynomials are defined as functions on the complex torus $\CC^n$.

\begin{theor}[(Kouchnirenko--Bernstein, \cite{bernst})]\label{thbernst} For every collection of finite sets $A=(A_1,\ldots,A_n)$ in $\Z^n$, there exists a proper exceptional algebraic set $B_A\subset \C^A=\C^{A_1}\oplus\ldots\oplus\C^{A_n}$, such that the number of common roots $x\in\CC^n$ of a system of polynomial equations $f_1(x)=\ldots=f_n(x)=0$ for every tuple of polynomials $(f_1,\ldots,f_n)\in\C^A$ outside $B_A$ equals the lattice mixed volume of (the convex hulls of)
$A_1,\ldots,A_n$.
\end{theor}

In the setting of the Kouchnirenko--Bernstein theorem, denote the mixed volume by $V$, then every loop in $\C^A\setminus B_A$, pointed at some tuple $f=(f_1,\ldots,f_n)$, defines a permutation of the roots of $f=0$. For all loops in $\C^A\setminus B_A$, these permutations form a subgroup of the group $S_V$ of all permutations of the $V$ roots of $f=0$. This subgroup will be called {\it the monodromy group of the general system of polynomial equations 
supported at} $A$ and denoted by $G_A$.

We shall be interested in the following two problems: 

\vspace{2ex}

(I) Compute $G_A$;

\vspace{2ex}

(II) Classify {\it solvable tuples} $A$, for which the multivalued function $\C^A\setminus B_A\to\CC^n$, assigning the roots of the system $f=0$ to an element $f\in\C^A\setminus B_A$, can be expressed by radicals.

\vspace{2ex}

The first problem helps to solve the second one, because a solvable tuple $A$ has a solvable monodromy group $G_A$ (see e.g. \cite{khtop}).

\begin{exa} For $n=1$ and $A=A_1=\{0,1,\ldots,d\}$, the problems above ask (I) for the monodromy of the generic univariate polynomial $c_dx^d+c_{d-1}x^{d-1}+\ldots+c_0$ and (II) for expressing its roots by radicals in terms of the coefficients $c_0, c_1,\ldots,c_d$. It is classically known that the monodromy $G_A$ equals $S_d$, and thus the general equation of degree $d$ is solvable for $d\leqslant 4$.
\end{exa}

For arbitrary $n$, the second problem, although not the first one, can be reduced without loss of generality to {\it reduced irreducible} tuples $A=(A_1,\ldots,A_n)$ in the sense of the following Definition \ref{Dred}. Thus, the subsequent Theorem \ref{thmain} leads to a complete solution of the problem (II), and seems to be a natural first step towards the solution of the problem (I).
\begin{defin}\label{Dred} 1. 
A tuple of finite sets $A_1,\ldots,A_k$ in $\Z^n$ is said to be reduced, if they cannot be shifted to the same proper sublattice of $\Z^n$. 

2. A tuple of finite sets $A_1,\ldots,A_k$ in $\Z^n$ is said to be irreducible  (resp. linearly independent), if it is impossible to shift all but $m$ of them (resp. $m-1$) to the same codimension $m$ sublattice for $m>0$.
\end{defin}
\begin{rem} 1. Mind the difference between reduced and reducible (i.e. non-irreducible).

2. Similar conditions were introduced by various authors for particular values of $n-k$ (c.f. for instance essential tuples in \cite{st94} for $k=n+1$). We prefer the names ``linearly independent'', ``reduced'' and ``irreducible'' (introduced in \cite{kh75} and \cite{abelruff} for $k=n$), because discriminants and other geometric objects, related to the system of equations $f=0$ for the general tuple $f\in\C^A$, tend to be reduced and irreducible in the sense of algebraic geometry if the tuple $A=(A_1,\ldots,A_k)$ has the property of the same name. See Remark \ref{rem1irr} and Theorem \ref{lirrdef0} for some instances of this correspondence.
\end{rem}
\begin{theor} \label{thmain} If $A=(A_1,\ldots,A_n)$ is a reduced irreducible tuple, 
then the monodromy group $G_A$ equals the symmetric group $S_V$.
\end{theor}
The proof is given at the end of this section. 

\vspace{2ex}

{\bf Systems of equations, solvable by radicals.} Since $S_V$ is not solvable for $V>4$, the preceding theorem implies the following.
\begin{sledst}[(Conjecture 1 in \cite{abelruff})] \label{Csolv} For a reduced irreducible tuple $(A_1,\ldots,A_n)$, the general system of equations supported at $(A_1,\ldots,A_n)$ is solvable by radicals if and only if it has at most 4 solutions, i.e. the lattice mixed volume of 
$A_1,\ldots,A_n$ does not exceed 4.
\end{sledst}
This fact actually gives the inductive classification of all solvable tuples $A=(A_1,\ldots,A_n)$. 

\vspace{2ex}

\begin{classif}\label{classif0} (0) We can and will assume without loss of generality that every $A_i$ contains 0. Indeed, otherwise shift $A_i$ by a vector $-a_i,\, a_i\in A_i$, to a set $\widetilde A_i$ containing $0$. Now, instead of polynomials $f_i\in\C^{A_i}$, we can study polynomials $f_i(x)/x^{a_i}\in\C^{\widetilde A_i}$, because they have the same roots as $f_i$. 

(1) We can and will assume that $A$ is reduced. Indeed, otherwise $A_i$ is the image of $B_i$ under a lattice embedding $j:\Z^n\to\Z^n$ for a reduced tuple $B=(B_1,\ldots,B_n)$, and we have the following fact: {\it the solvability of $B$ is equivalent to the solvability of $A$.} 

{\sc Proof.} Consider the surjection of complex tori $h:\CC^n\to\CC^n$, corresponding to the embedding $j$ of their character lattices, so that $h(x)^b=x^{j(b)}$ for $x\in\CC^n$ and $b\in\Z^n$, then every tuple of polynomials $f\in\C^A$ has the form $f(x)=g(h(x)),\, g\in\C^B$. Since $h$ is 
invertible by radicals, then $f=0$ and $g=0$ are solvable by radicals simultaneously. $\quad\square$

(2) We can and will assume that $A$ is irreducible. Otherwise, up to reordering, the sets $A_1,\ldots,A_k,\, 0<k<n$, belong to the same $k$-dimensional plane $L\subset\Z^n$, and, denoting the tuple of the images of the other $A_i$'s under the projection $\Z^n\to\Z^n/L$ by $A''$, we have the following fact: {\it the solvability of $A$ is equivalent to the solvability of the smaller dimensional tuples  $A'=(A_1,\ldots,A_k)$ and $A''$}. 

{\sc Proof.} Note that upon an appropriate automorphism of $\CC^n$, the polynomial $f_i\in\C^{A_i}$ depends only on the first $k$ coordinates for $i\leqslant k$, so, substituting these coordinates with a solution of $f_1=\ldots=f_k=0$ in the system of equations $f_{k+1}=\ldots=f_n=0$, we obtain a system of the form $g=0,\, g\in\C^{A''}$. Thus solving a generic system $f=0$ supported in $A$ amounts to solving a generic system $f_1=\ldots=f_k=0$ supported in $A'$ and a system $g=0$, which is also generic in $\C^{A''}$ in the sense that assigning $g$ to $f$ is a  dominant map $\C^A\to\C^{A''}$. $\quad\square$

(3) Finally, a reduced and irreducible tuple $A$ is solvable if and only if the lattice mixed volume of 
$A_1,\ldots,A_n$ does not exceed 4 (by Corollary \ref{Csolv}).
\end{classif}

This algorithm reduces the classification of solvable systems of equations to the classification of irreducible mixed volume 4 tuples of lattice sets. The latter classification is given in \cite{abelruff} in dimension 2, and is moreover finite in every dimension, see Theorem \ref{thmv} below for details.
\begin{rem} In the same way, the classification of systems of equations solvable by $k$-radicals in the sense of \cite{khtop} (i.e. the ones that can be reduced to solving univariate polynomial equations of degree at most $k$) is reduced to the classification of tuples of lattice sets of mixed volume at most $k$.
\end{rem}
\begin{exa} \label{exa11} For $n=2$, if a reduced consistent general system of equations is solvable by radicals, then its Newton polygons either have lattice mixed volume at most four (there are 14 such maximal pairs up to automorphisms of $\Z^2$, see \cite{abelruff}), or equal a segment $I$ of lattice length at most 4 and an arbitrary polygon $P$, whose support lines parallel to $I$ are at the lattice distance not exceeding 4 from each other.
\end{exa}

{\bf Classification of small polytopes.} Each of the infinitely many pairs $(I,P)$ in the preceding example has mixed volume at most 16, due to the following fact. We denote the lattice mixed volume of the convex hulls of $A_1,\ldots,A_n$ by $\MV(A_1,\ldots,A_n)$.
\begin{theor} \label{thmvpr} Let $B_1,\ldots,B_N$ be lattice sets in $\Z^N$ and $A_1,\ldots,A_M$ in $\Z^N\oplus\Z^M$. 
Then $\MV(A_1,\ldots,A_M,B_1,\ldots,B_N)=\MV(p A_1,\ldots,p A_M)\MV(B_1,\ldots,B_N)$, where $p:\Z^N\oplus\Z^M\to\Z^M$ is the standard projection.
\end{theor}
This well known fact admits an especially simple proof in the spirit of Classification \ref{classif0}.2 (\cite{sb06}, see e.g. Lemma 4 in \cite{theo} for a geometric proof).

\noindent {\sc An algebraic proof.} 
For $f_i\in\C^{A_i}$ and $g_j\in\C^B_j$, every solution of the system $f=g=0$ is of the form $(x_0,y_0)$, where $x_0\in\CC^N$ is a solution of the system $g=0$ and $y_0\in\CC^M$ is a solution of the system $f(x_0,y)=0$. For generic $f$ and $g$, the number of solutions of the three mentioned systems equals the three lattice mixed volumes in the statement by the Kouchnirenko--Bernstein theorem.
$\quad\square$

\vspace{2ex}

This reduces the infinite classification of tuples with small mixed volume to the classification of irreducible tuples, which is already finite.
\begin{theor}\label{thmv} For every $n$ and $V$, there are finitely many irreducible tuples $(A_1,\ldots,A_n)$ in $\Z^n$ of mixed volume $V$, up to automorphisms of $\Z^n$ and shifts of the sets.
\end{theor}
The proof is given in Section \ref{Smv}. Moreover, if we restrict our attention to the {\it unmixed} case, where $A_1=\ldots=A_n=A$, the classification becomes essentially finite across all dimensions: it was shown in \cite{abelruff} that every reduced $A\subset\Z^n$ of lattice volume 4 can be obtained from 34 ``elementary'' configurations of dimension at most 6  by affine automorphisms of $\Z^n$ and constructing cones over lattice sets in the following sense.
\begin{defin} The cone over $B\subset\Z^m$ is the set $c(B)=\{0\}\cup \bigl(B\times\{k\}\bigr)\subset\Z^{m+1}$.
\end{defin}
\begin{rem} 1. The same is true for every value of the volume, as shown in \cite{n16}, Corollary 3.1 (although, starting from volume 5, the classification of non-cones seems to be incomprehensibly large).

2. In the notation of the preceding definition, the solution by radicals of the system $f_0=\ldots=f_m=0$ supported at the cone $c(B)$ can be reduced to the solution by radicals of the system $g_1=\ldots=g_m=0$ supported at its base $B$, by setting $g_i(x)=f_i(x)/f_{i,0}-f_0(x)/f_{0,0}$, where $f_{i,0}$ is the constant term of $f_i$. Thus the solution by radicals of all solvable unmixed systems of arbitrarily many variables reduces to the 34 elementary ones, listed in \cite{abelruff}.

3. The classification of the 34 non-cones of volume four in \cite{abelruff} includes only reduced ones (or spanning ones, in terms of \cite{n16}), because this suits the needs of Corollary \ref{Csolv}. The classification of all (possibly non-reduced) non-cones of volume four is also possible, but is more complicated and not finite due to empty simplices, see \cite{v4}.
\end{rem}

\vspace{2ex}

{\bf Monodromy of reducible systems of equations.} In contrast to the problem of solvability, the computation of the monodromy of an arbitrary tuple cannot be reduced to the case of reduced irreducible tuples easily. We formulate a conjecture regarding non-reduced tuples and show by an example that the case of reducible tuples is yet more complicated (so that we do not even make any predictions).

\begin{conj}\label{conjwreath} In the setting of Step (1) of Classification \ref{classif0}, if the tuple $B$ is reduced and irreducible of mixed volume $d$, then the monodromy group $G_A$ equals the wreath product of $\coker j$ and $S_d$ acting on $\{1,\ldots,d\}$.
\end{conj}
\begin{rem}\label{remconj} We now explain why $G_A$ obviously embeds into this wreath product, so the problem is whether the embedding is actually an isomorphism. In the notation of Part 1 of Classification \ref{classif0}, the roots of $f=0$ split into the fibers of the surjection $h:\{f=0\}\to\{g=0\}$. All fibers are cosets of the subgroup $\coker j\subset\CC^n$, and every monodromy permutation of the set $\{f=0\}$ ``respects $j$'', i.e. it sends every fiber into a fiber, preserving its $(\coker j)$-torsor structure. In particular, the group $G_A$ is contained in the group of all permutations respecting $j$, and the latter is exactly the sought wreath product.
\end{rem}
\begin{exa}\label{excayley} If the tuple $A$ is as shown on the left (Figure 1), then $G_A$ is obviously equal to $V_4\subset S_4$, generated by $(12)(34)$ and $(13)(24)$. However, if the tuple $A$ is as shown on the right, then its Cayley discriminant (Definition \ref{defdiscr3}) has codimension 1, so a small loop around this discriminant corresponds to a transposition in $G_A$ (see Remark \ref{forexa1}), thus the group is strictly greater than $V_4$ (actually, it equals $D_8$). This is despite, in the notation of Step (2) of Classification \ref{classif0}, the groups $G_{A'}$ and $G_{A''}$ are the same (equal to $S_2$) for both examples. Thus $G_A$ is not defined solely by $G_{A'}$ and $G_{A''}$.

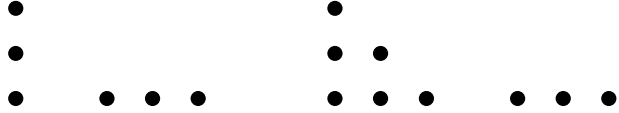
\begin{figure}
\centering
\psscalebox{0.5 0.5} 
{
\begin{pspicture}(0,-1.3971155)(15.994231,1.3971155)
\psdots[linecolor=black, dotsize=0.4](0.19711533,1.2)
\psdots[linecolor=black, dotsize=0.4](0.19711533,0.0)
\psdots[linecolor=black, dotsize=0.4](0.19711533,-1.1999999)
\psdots[linecolor=black, dotsize=0.4](2.5971153,-1.1999999)
\psdots[linecolor=black, dotsize=0.4](3.7971153,-1.1999999)
\psdots[linecolor=black, dotsize=0.4](4.997115,-1.1999999)
\psdots[linecolor=black, dotsize=0.4](8.5971155,-1.1999999)
\psdots[linecolor=black, dotsize=0.4](8.5971155,0.0)
\psdots[linecolor=black, dotsize=0.4](8.5971155,1.2)
\psdots[linecolor=black, dotsize=0.4](9.797115,0.0)
\psdots[linecolor=black, dotsize=0.4](9.797115,-1.1999999)
\psdots[linecolor=black, dotsize=0.4](10.997115,-1.1999999)
\psdots[linecolor=black, dotsize=0.4](13.397116,-1.1999999)
\psdots[linecolor=black, dotsize=0.4](14.5971155,-1.1999999)
\psdots[linecolor=black, dotsize=0.4](15.797115,-1.1999999)
\end{pspicture}
}
\caption{Two reducible tuples}
\end{figure}
\end{exa}

Nevertheless, we can confirm in our setting the ``symmetric or imprimitive'' dichotomy, conjectured in \cite{schsot} for Schubert enumerative problems, modulo one obvious exclusion.
\begin{exa} \label{exaexcl} Let $B$ and $C$ be tuples of finite sets of lattice mixed volume 1 in $\Z^k$ and $\Z^m$ respectively, $k>0,\, m\geqslant 0$ (see \cite{EG1} or Section 2 below for the classification of such tuples), and let $j:\Z^k\to\Z^k\oplus\Z^m$ send $v$ to $(pv,0)$ for some odd prime $p$. Let $P'$ be the tuple $j(B)$ in $\Z^k\oplus\Z^m$, and let $P''$ be a tuple of $m$ sets in $\Z^k\oplus\Z^m$ whose projections to $\Z^m$ form the tuple $C$. Then the mixed volume of a tuple $P=(P',P'')$ equals $p$, and, moreover, by Remark \ref{remconj}, the monodromy group $G_P$ is a subgroup of $\Z/p\Z$, i.e. equals $\Z/p\Z$ or the trivial group, of which the former is primitive and the latter is not. Actually one can check that $G_P$ always equals $\Z/p\Z$ in accordance with Conjecture \ref{conjwreath} (which is obvious in the 1-dimensional case, i.e. for the equation $c_px^p+c_0=0$, corresponding to $P_1=\{0,p\}\subset\Z^1$, and less obvious in general).

A tuple that can be identified with $P$ by an isomorphism of lattices will be called a prime tuple.
\end{exa}
\begin{defin} \label{defnum} A tuple of sets $A=(A_1,\ldots,A_n)$ in $\Z^n$ is said to be numerically non-reduced, if there exist sets $B_1,\ldots,B_k$ in $\Z^k$ and an embedding $j:\Z^k\to\Z^n$, such that the lattice mixed volume of $B_1,\ldots,B_k$ is greater than 1, the embedding is not saturated (i.e. $\Z^n/j(\Z^k)$ is not free), and $j(B_1),\ldots,j(B_k)$ coincide with $k$ of the sets $A_1,\ldots,A_n$ up to a shift.

The tuple $A$ is said to be numerically reducible, if $k<n$ of $A_i$'s can be shifted to a $k$-dimensional sublattice $L$ such that the lattice mixed volumes of  both $A'=($the tuple of $A_i$'s shifted to $L)$ and $A''=($the tuple of the images of the rest of $A_i$'s under the projection $\Z^n\to\Z^n/L)$ are greater than 1.
\end{defin}
The name is chosen because the mixed volume $V$ of the tuple $A$ equals the product of the mixed volumes of $A'$ and $A''$ by Theorem \ref{thmvpr}.
\begin{theor} For every non-prime tuple $A$ (see Example \ref{exaexcl}), the monodromy group $G_A$ is the symmetric group $S_V$ 
if the tuple $A$ is numerically reduced and irreducible, and is imprimitive otherwise.
\end{theor}
{\sc Proof.} If the tuple $A$ of subsets of $\Z^n$ is numerically non-reduced, then, in the notation of Definition \ref{defnum}, let $h:\CC^n\to\CC^k$ be the surjection of tori, corresponding to the embedding $j:\Z^k\to\Z^n$ of their character lattices so that $h(x)^b=x^{j(b)}$ for all $x\in\CC^n$ and $b\in\Z^k$. Then every system of equations $f(x)=0,\, f\in\C^A$ contains a subsystem of the form $g(h(x))=0,\, g\in\C^B$. By Remark \ref{remconj}, the fibers of the surjection $h:\{f=0\}\to\{g=0\}$ are blocks of the monodromy action of $G_A$. The number and size of the blocks are greater than 1, because the mixed volume of $B_1,\ldots,B_k$ is greater than 1, and $j$ is not saturated.

If $A$ is numerically reducible, then, in the notation of Definition \ref{defnum}, upon an appropriate automorphism of $\CC^n$ and reordering the tuple, we may assume that $A_1,\ldots,A_k$ are contained in the first $k$-dimensional coordinate plane $L\subset\Z^n,\, 0<k<n$, and the mixed volumes $V'$ and $V''$ of both $A'=(A_1,\ldots,A_k)$ and $A''=($the images of $A_{k+1},\ldots,A_n$ in $\Z^n/L)$ are greater than 1. In this case, every common root of a generic tuple of polynomials $f=(f_1,\ldots,f_n)\in\C^A$ is of the form $(x',x'')$, where $x'\in\CC^k$ is one of the $V'$ roots of the system $f'=(f_1,\ldots,f_k)$. In particular, the fibers of the projection $\{f=0\}\to\{f'=0\}$ are $V'>1$ blocks of size $V''>1$ for the action of the monodromy group $G_A$, so this action is imprimitive.

If the tuple $A$ is numerically reduced, numerically irreducible and not prime, then it is reduced. So, if $A$ is irreducible in this case, then $G_A$ is symmetric by Theorem \ref{thmain}.

Thus, it remains to consider reducible $A$ that is numerically reduced, numerically irreducible and not prime. In this case, in the notation of Part 2 of Classification \ref{classif0}, the tuples $A'$ and $A''$ are also numerically reduced, numerically irreducible and not prime, and the mixed volume of one of them equals $1$. Thus $G_A$ equals the monodromy group of the other one, which is symmetric by induction on the dimension.
$\quad\square$

\vspace{2ex}

{\bf Structure of the paper.} In Section \ref{Smv}, we prove and discuss Theorem \ref{thmv}. The rest of the paper is devoted to the proof of Theorem \ref{thmain}. In Section \ref{Sdiscr}, we reduce the assumption of irreducibility to a more general notion of dual effectiveness (antonym to dual defectiveness, see Definition \ref{defdd} below). 
\begin{theor} \label{thirrdef} A reduced irreducible tuple of $n$ sets in $\Z^n$ is dual effective unless, upon an automorphism of the lattice, 
all of its sets can be shifted to the standard simplex (i.e. the system of equations is essentially linear).
\end{theor}
For the proof, see Corollary \ref{lirrdef}. Besides the relation to Galois theory, this result may be important as an illustration of a new approach to dual defectiveness in the toric setting, independent of the known ones \cite{dr}, \cite{dfs}, \cite{cc}, \cite{jems}, \cite{fi}, \cite{fors}.
\begin{rem} 1. In the case of full-dimensional tuples, Theorem \ref{thirrdef} was deduced from \cite{fi} in \cite{bn}, settling the conjecture from \cite{cdd}. Our proof is independent of \cite{fi}, and it would be important to extend the technique of \cite{bn} from full-dimensional tuples to irreducible ones.

2. It would be important to drop the irreducibility assumption and completely classify dual defective tuples in various senses (see Remark \ref{qred}), as Example \ref{excayley} suggests.
\end{rem}

\begin{theor} \label{thtransp} If $A$ is a reduced dual effective tuple, then the monodromy $G_A$ contains a transposition.
\end{theor} 
Roughly speaking, the transposition is produced by running a small loop around the discriminant, see Theorem \ref{luseful} for the proof and Theorem \ref{luseful2} for a possible generalization to non-square systems of equations.
\begin{theor} \label{th2tr} 
If $A$ is a reduced irreducible tuple, then the monodromy $G_A$ is doubly transitive.
\end{theor}
The proof is standard and is given in Section \ref{S2tr}. 

\vspace{2ex}

{\sc Proof of Theorem \ref{thmain}.} 
Unless the system of equations generically has one solution (satisfying $G_A=S_1$), Theorem \ref{thirrdef} ensures that the tuple is dual effective, so the monodromy contains a transposition by Theorem \ref{thtransp}. Since it is also doubly transitive by Theorem \ref{th2tr}, it coincides with the symmetric group. $\quad\square$

\section{Lattice polytopes of small mixed volume} \label{Smv}
\begin{theor}[\cite{LZ}]\label{thmain0} For any $n$, there are finitely many convex lattice polytopes of a given lattice volume in $\Z^n$ up to affine automorphisms of the lattice.
\end{theor}

\begin{theor}[(Minkowski, \cite{Mink})]\label{thmink} A tuple is linearly dependent if and only if its mixed volume equals 0.
\end{theor}
{\sc Proof of Theorem \ref{thmv}.} Tuples $(B_1,B_1,B_3,\ldots,B_N)$ and $(B_2,B_2,B_3,\ldots,B_N)$ are said to be AF-descendants of $(B_1,B_2,B_3,\ldots,B_N)$, if both of them are linearly independent. If the tuple $B'$ is the AF-descendant of $B$, then, by the Aleksandrov--Fenchel inequality and Theorem \ref{thmink}, we have $$\MV B'\leqslant (\MV B)^2.\eqno{(*)}$$

Every linearly independent tuple $B$ that entirely consists of sets contained in the irreducible tuple $A$, can be obtained from $A$ by taking a sequence of AF-descendants $A',A'',\ldots,A^{(k)}=B$.\footnote[4]{July 2020: Apparently, I cannot prove this statement (the existence of  a sequence of AF-descendants) at the request of the readers, so I include an addendum at the end of the text to circumvent this step of the proof. This patch does not change the statement of Theorem 1.11, other steps of its proof and other parts of the paper.} Applying the inequality $(*)$ to this sequence, we conclude: if all sets of the tuple $B$ are contained in the irreducible tuple $A$, then $$\MV B\leqslant (\MV A)^{2^N}. \eqno{(**)}$$
Note that $(**)$ trivially holds also for linearly dependent tuples $B$ by Theorem \ref{thmink}.

We can now estimate the lattice volume of the Minkowski sum $A_1+\dotsb+A_N$ as follows: write it as $\MV(A_1+\dotsb+A_N,\ldots,A_1+\dotsb+A_N)$, open the brackets and estimate every term by the inequality $(**)$. 
As a result, for every irreducible tuple $(A_1,\ldots,A_N)$ of mixed volume $V$, the volume of the Minkowski sum $A_1+\dotsb+A_N$ is at most $N^N V^{2^N}$, so by Theorem \ref{thmain0} there are finitely many possibilities for $A_1+\dotsb+A_N$ and hence for $(A_1,\ldots,A_N)$.
$\quad\square$
\begin{rem} It would be interesting to obtain a sharper estimate on the volume of $A_1+\ldots+A_N$ in terms of the mixed volume of an irreducible tuple $(A_1,\ldots,A_N)$.
\end{rem}
The classification of irreducible tuples is known only up to mixed volume 4 in dimension 2 (see \cite{abelruff}), and up to mixed volume 1 in arbitrary dimension:
\begin{sledst}[(Minkowski)] The unique irreducible tuple of mixed volume 0 is a point in $\Z^1$.
\end{sledst}
\begin{theor}[\cite{EG1}]\label{thmv1} The unique (up to automorphisms of the lattice and shifts of polytopes) maximal (by inclusion) irreducible tuple of lattice polytopes of mixed volume 1 in $\Z^N$ is the tuple of $N$ copies of the standard simplex.
\end{theor}

\section{Discriminants and dual defectiveness} \label{Sdiscr}

$\quad$

{\bf Mixed resultants.} Let $A=(A_0,\ldots,A_n)$ be a tuple of finite sets in $\Z^n$.
\begin{defin} The $A$-resultant $R_A$ is the closure of the set of all tuples of polynomials $f=(f_0,\ldots,f_n)\in\C^A$ that have a common root $f_0(x)=\ldots=f_n(x)=0,\, x\in\CC^n$.
\end{defin}
\begin{exa} For $n=1$, the set $R_A$ is the zero locus of the classical Sylvester resultant.
\end{exa}
\begin{theor}[(Theorem 2.26 in \cite{stekl}, see also \cite{st94} for the first part of the statement)] \label{thres}
If $A$ is irreducible, then the resultant $R_A$ is a nonempty irreducible hypersurface, and a generic tuple $f\in R_A$ has a unique common root in $\CC^n$.
\end{theor}

{\bf Gelfand--Kapranov--Zelevinsky discriminants.} Let $A\subset\Z^n$ be a finite set.
\begin{defin}[\cite{gkz}] The $A$-discriminant $D_A$ is the closure of the set of all polynomials $f\in\C^A$ that have a singular root $f(x)=0,\, df(x)=0,\, x\in\CC^n$.
\end{defin}
\begin{exa} For $n=1$, the set $D_A$ is the zero locus of the classical discriminant.
\end{exa}
\begin{defin}\label{defdualdef1} The tuple $A$ is said to be dual defective if $D_A$ is not a hypersurface, and dual effective otherwise.
\end{defin}
This is equivalent to the projectively dual variety to the toric variety $X_A$ is not a hypersurface (hence the name). The study of dual defective projective varieties is a classical topic in algebraic geometry \cite{Ein85}. In particular, there is an extensive literature on the classification of dual defective lattice sets, see \cite{dr}, \cite{jems} and \cite{fi} for some of the most explicit answers (the first one is for the case of smooth toric varieties). 
\begin{exa} The set $A=\{(00), (10), (20), (01)\}\subset\Z^2$ is defective.
\end{exa}
\begin{theor}[\cite{gkz}]\label{thgkz} If a dual effective $A$ cannot be shifted to a proper sublattice of $\Z^n$, then a generic polynomial $f\in D_A$ has a unique singular root $x\in\CC^n$, and the Hessian of $f$ at this root is non-degenerate. 
\end{theor} 
If $A$ is dual effective, then the set $D_A$ is the zero locus of a unique irreducible integer polynomial on $\C^A$ (up to the choice of the sign). This polynomial is also called the $A$-discriminant. The coefficients $c_a, a\in A$, of the general Laurent polynomial $\sum_{a\in A} c_ax^a$ in $\C^A$ form the natural system of coordinates in $\C^A$, and we shall consider the $A$-discriminant as a polynomial of $c_a, a\in A$.
\begin{lemma}[(Lemma 2.21 in \cite{dcg})] For every dual effective $A\subset\Z^n$ and every $a\in A$, the $A$-discriminant has positive degree in $c_a$.
\end{lemma}
\begin{rem}\label{remforget} For every $B\subset A$, there is a natural forgetful projection $\C^A\to\C^B$, sending $\sum_{a\in A}c_ax^a$ to $\sum_{a\in B}c_ax^a$, and we shall denote the preimage of $D_B$ under this map also by $D_B$.
\end{rem}
\begin{sledst}\label{cordiff} If $A$ is dual effective, then $D_A\ne D_B$ for every $B\subsetneq A$.
\end{sledst}

{\bf Discriminants of systems of equations.} For a tuple $A=(A_1,\ldots,A_k)$ of finite sets in $\Z^n,\, 2\leqslant k\leqslant n$, the concept of the discriminant of the system of equations supported at $A$ is ambiguous. We introduce three different versions of this notion that appear in the literature, and it will be important for us that all of them coincide for irreducible tuples. Denote the standard basis in $\Z^k$ by $e_1,\ldots,e_k$, and, for every $I\subset\{1,\ldots,k\}$, let $A_I$ be the {\it Cayley configuration} $\cup_{i\in I} \{e_i\}\times A_i\subset\Z^k\times\Z^n$. For every $f\in\Z^A$, let $f_I$ be the polynomial $\sum_{i\in I}\lambda_i f_i(x)\in\C^{A_I}$ of variables $\lambda=(\lambda_1,\ldots,\lambda_k)\in\CC^k$ and $x\in\CC^n$.
\begin{defin} \label{defdiscr3}
1) The {\it naive $A$-discriminant} \cite{dcg} is the closure of the set of all tuples $f\in\C^A$ having a singular common root $x\in\CC^n$ (so that $f_1(x)=\ldots=f_k(x)=0$ and $df_1(x),\ldots,df_k(x)$ are linearly dependent).

2) The {\it mixed $A$-discriminant} (\cite{cdd} for $k=n$) is the closure of the set of all tuples $f\in\C^A$ having a {\it non-degenerate singular common root} $x\in\CC^n$ (i.e. a singular common root such that no proper subtuple of $df_1(x),\ldots,df_k(x)$ is linearly dependent).

3) The {\it Cayley $A$-discriminant} \cite{dcg} is the image of the discriminant $D_{A_{\{1,\ldots,k\}}}\subset\C^{A_{\{1,\ldots,k\}}}$ under the natural isomorphism $\C^{A_{\{1,\ldots,k\}}}\to\C^A$ inverse to sending every $f$ to $f_{\{1,\ldots,k\}}$.
\end{defin}
All of these sets obviously coincide for the Gelfand--Kapranov--Zelevinsky case $k=1$. However, for $k>1$ (including $k=n$), they may be pairwise different (Example 1.2 in \cite{cdd}), and have irreducible components of different dimensions (Example 2.25 in \cite{dcg}). 
Nevertheless, this difference disappears for irreducible tuples.
\begin{theor} \label{thcoinc} If 
$A$ is irreducible, the three discriminant sets of Definition \ref{defdiscr3} coincide up to irreducible components of codimension greater than 1.
\end{theor}
{\sc Proof.} 
If the Cayley discriminant has codimension greater than 1, then so does the naive discriminant by Theorem 2.31 in \cite{dcg} and the mixed discriminant (as its subset). 

To study the opposite case, define $\Sigma_{\{j_1,\ldots,j_q\}}$ as the set of all tuples $f=(f_1,\ldots,f_k)\in\C^A$ such that $f(x)=0$ for some $x\in\CC^n$ and $\sum_i \lambda_idf_{q_i}(x)=0$ for some $(\lambda_1,\ldots,\lambda_q)\in\CC^q$ (Definition 2.33 in \cite{dcg}).

If the Cayley discriminant has codimension 1, then this hypersurface $H$ is the only codimension 1 component of the naive discriminant  by Theorem 2.31 in \cite{dcg} and the only codimension 1 set of the form $\Sigma_J$ (namely, the one corresponding to $J=\{1,\ldots,k\}$) by Lemma 2.34 in \cite{dcg}. The latter fact implies that a singular common root of a generic tuple $f\in H$ is non-degenerate (because the linear dependence of its differentials $df_j$ for $j\in J'\ne\{1,\ldots,k\}$ would imply that $\Sigma_{J'}=\Sigma_J$ also has codimension 1).
Thus $H$ is also a codimension 1 component of the mixed discriminant, and the latter has no other codimension 1 components, because it is contained in the naive discriminant.
$\quad\square$

\vspace{1ex}

{\bf Dual defectiveness of systems of equations.}
\begin{defin} \label{defdd} By Theorem \ref{thcoinc},  for an irreducible tuple $A$, we can denote the common hypersurface components of the three discriminant sets of Definition \ref{defdiscr3} by $D_A$, and call this hypersurface the $A$-discriminant. The irreducible tuple $A$ is said to be dual defective if $D_A$ is empty, and dual effective otherwise.
\end{defin}
If the tuple $A$ consists of one set $A_1\subset\Z^n$, then it is irreducible, and its dual defectiveness is the same property as in Definition \ref{defdualdef1}.
\begin{conj} For irreducible tuples, the three discriminant sets of Definition \ref{defdiscr3} coincide completely, i.e. they are the same irreducible set. 
\end{conj}
\begin{rem}\label{qred} The interrelation between the three notions of the discriminant in Definition \ref{defdiscr3} is not completely understood for reducible tuples so far, particularly concerning the higher codimension components. As a consequence, the notion of dual defectiveness for reducible tuples splits into several non-equivalent versions, looking for the non-existence of codimension 1 components and/or existence of higher codimension components in any of the three notions of the discriminant. It would be important to understand how these numerous versions are related.
\end{rem}
\begin{rem}\label{rem1irr} As we have observed, the irreducibility of the tuple $A$ implies the irreducibility of the codimension 1 part of the naive $A$-discriminant. On the other hand, the codimension 1 part of the naive $A$-discriminant tend to be reducible if $A$ is reducible, see Lemma 2.34 in \cite{dcg}. The situation with reduced tuples is similar: the codimension 1 components $D_i$ of the naive $A$-discriminant come with natural multiplicities equal to the number of singular roots of the system $f=0$ for a generic tuple $f\in D_i$. By Theorem \ref{thgkz}, an irreducible tuple $A$ is reduced if an only if $D_A$ is reduced in the sense of the aforementioned multiplicity (see \cite{eadv} for the computation of the multiplicities for non-reduced and reducible tuples).
\end{rem}
\begin{lemma}\label{lisol} An irreducible tuple $A$ is dual effective if and only if some $f\in\C^A$ has an isolated singular root.
\end{lemma}
{\sc Proof.} If the tuple $f=(f_1,\ldots,f_k)$ has an isolated singular root $x$, then the set of tuples in $\C^A$ that have a singular root contains a hypersurface in a small neighborhood of $f$. Indeed, the  projection $\pi$ of the incidence set $\{(\tilde x,\tilde f)\,|\, \tilde f(\tilde x)=0\}\subset\CC^n\times\C^A$ to $\C^A$ has the critical set $C$ of dimension by 1 smaller than $\C^A$. Since $x$ is an isolated singular root of $f$, the fibers of the projection $\pi:C\to\C^A$ near $(x,f)$ are finite, thus the image of $C$ contains a hypersurface passing through $f$. Thus, according to the naive version of the definition of the discriminant (see Definition \ref{defdiscr3}), $D_A$ contains a non-empty hypersurface.

To prove the statement in the other direction, recall that $f_{\{1,\ldots,k\}}$ is a homogeneous polynomial in the variables $\lambda_1,\ldots,\lambda_k$, so the equation $f_{\{1,\ldots,k\}}=0$ defines a subset in $\CP^{k-1}\times\CC^n$. Denote the image of the torus $\CC^k$ under the projection $\C^k\to\CP^{k-1}$ by $T$. 

In this notation, if the tuple $A$ is dual effective, then so is $A_{\{1,\ldots,k\}}$, then, by Theorem \ref{thgkz}, a generic polynomial $f_{\{1,\ldots,k\}}$ in it has a unique (and thus isolated) singular root in $T\times\CC^n$, then so does the tuple $f\in\C^A$. $\quad\square$

\vspace{1ex} 

{\bf The proof of Theorem \ref{thirrdef}.} For a finite set $A\subset\Z^n$, let $\CP^{A}$ be the projective space with the homogeneous coordinates $z_a,\,\ a\in A$, and let $m=m_{A}:\CC^n\to\CP^{A}$ be the monomial map such that $m_{A}(x)$ has coordinates  $z_a=x^a$.
\begin{defin} The $A$-image of an algebraic set $V\subset\CC^n$ is the image of $m_A(V)$ in $\CP^{A}$.
\end{defin}
\begin{rem}
The $A$-image is usually not closed. In what follows, whenever we discuss its degree and irreducibility, we refer to the corresponding properties of its closure. On the other hand, its projectively dual set is defined as the set of all tangent hyperplanes to its smooth points, and is usually also not closed.\end{rem}
\begin{theor}\label{lirrdef0} Let  $A=(A_1,\ldots,A_k)$ and $A'=(A_2,\ldots,A_k)$ be tuples of finite sets in $\Z^n, n\geqslant k$, and let $M$ be the $A_1$-image of the complete intersection $f=0$ for a generic tuple of polynomials $f=(f_2,\ldots,f_k)\in\C^{A'}$.

1) If $A$ is irreducible, then $f=0$ and $M$ are irreducible.

2) If $A$ is moreover reduced, then $M$ is also reduced (in the sense that the map $m=m_{A_1}:\{f=0\}\to M$ has degree 1).

3) Assume that $A$ is reduced and irreducible. Then the degree of $M$ is greater than 1 unless the sets $A_1,\ldots,A_k$ can be shifted to the same lattice simplex of lattice volume 1.

4) Assume that $A$ is reduced and irreducible. Then $A$ is dual defective if and only if $M$ is dual defective (i.e. its projectively dual set has codimension greater than 1). Moreover, if $A$ and $M$ are dual effective, then a generic tuple of polynomials in the discriminant $D_A$ has a unique singular root, and this root is non-degenerate.
\end{theor}
\begin{rem}\label{remmiln} 1. We shall apply this lemma for $k=n$, in which case by the non-degenerate singular root we mean just the root of multiplicity 2. However, Part 4 makes sense for arbitrary $k\leqslant n$. In this case a root $x$ of a system of equations $g=0$ is said to be singular non-degenerate, if $g=0$ defines an isolated singularity of a complete intersection in a neighborhood of $x$, and its Milnor number equals 1 (see \cite{looijenga}).

2. Part 1 for $k<n$ actually takes place and will be proved under a strictly weaker assumption that we call coirreducibility (cf. Definition \ref{Dred}): no $m$ sets of the tuple $A'$ can be shifted to the same $m$-dimensional sublattice.
\end{rem}
{\sc Proof.} If $A$ is irreducible and $k\leqslant n$, then the tuple $A'$ is coirreducible (in the sense of Remark \ref{remmiln}.2), so Part 1 follows from \cite{kh15} for $f=0$ and thus also for $M$. 

We shall now assume without loss of generality that $0\in A_1$, because all properties of $A$ mentioned in the statement are invariant under parallel translations. For every linear form $l$ on $\CP^{A_1}$, denote the rational function $l/z_0$ by $\tilde l$. In this notation, assigning the function $f_1(x)=\tilde l(m(x))$ to a form $l$ (or, in coordinates, assigning the polynomial $f_1(x)=\sum_{a\in A_1} c_ax^a$ to the form $l(z)=\sum_{a\in A_1} c_az_a$), we establish an isomorphism between $\C^{A_1}$ and the space of linear forms on $\CP^{A_1}$.

Assume towards contradiction that Part 2 does not hold. Then, for generic linear forms $l_1,\ldots,l_{n-k+2}$ such that the plane $l_\bullet=0$ intersects $M$, an intersection point would have more than one preimage in $f=0$, i.e. a generic tuple of polynomials $$\Big( \tilde l_1(m(\cdot)),\ldots, \tilde l_{n-k+2}(m(\cdot)), f_2,\ldots,f_k\Big) \in R_{B}$$ supported at the tuple $$B=(\underbrace{A_1,\ldots,A_1}_{n-k+2},A_2,A_3,\ldots,A_k)$$ would have more than one common root. This would contradict Theorem \ref{thres}, because irreducibility of $A$ implies irreducibility of $B$.

In the setting of Part 3, we may assume without loss of generality by Theorem \ref{thmv1} that the tuple $$B'=(\underbrace{A_1,\ldots,A_1}_{n-k+1}
,A_2,A_3,\ldots,A_k)$$ has mixed volume greater than 1, because it is reduced and irreducible. Then the degree of $M$ is greater than 1, because it equals the number of intersections of $V$ with a generic plane $l_1=\ldots=l_{n-k+1}=0$, i.e. the number of common roots of a generic tuple of polynomials $\big(\tilde l_1(m(\cdot)),\ldots,\tilde l_{n-k+1}(m(\cdot)),f_2,\ldots,f_k\big)\in\C^{B'}$, which equals the mixed volume of $B'$ by the Kouchnirenko--Bernstein theorem.

It remains to prove Part 4. If $M$ is dual effective, then the hyperplane $l=0$, corresponding to a smooth point of the projectively dual variety, is tangent to $M$ at a unique point $z$, and the tangency is non-degenerate (in the sense that the restriction of $\tilde l$ to $M$ has the non-degenerate Hessian at $z$). Then the restriction of the polynomial $f_1(x)=\tilde l(m(x))$ to the complete intersection $f=(f_2,\ldots,f_k)=0$ has a unique and non-degenerate singular root, then the resulting tuple $(f_1,f_2,\ldots,f_k)$ has a unique and non-degenerate singular root. By Lemma \ref{lisol}, this implies that $A$ is dual effective. The other direction is proved in the same way.
$\quad\square$

\begin{sledst}[(refined Theorem \ref{thirrdef})] \label{lirrdef} A reduced irreducible tuple of sets $A=(A_1,\ldots,A_n)$ in $\Z^n$ is dual effective unless, upon an automorphism of the lattice, 
all of its sets can be shifted to the standard simplex. Moreover, in this case a generic tuple $f\in D_A$ has a unique multiple root, and this root has multiplicity 2.
\end{sledst}
{\sc Proof.} By Theorem \ref{lirrdef0}.1-3, the closure of $M$ is a reduced irreducible curve of degree greater than 1. Since every such curve is dual effective, the sought statement follows from Theorem \ref{lirrdef0}.4.$\quad\square$
\begin{rem} Excluding the notion of the projectively dual variety from this reasoning, we can describe more explicitly the picture in $\CP^{A_1}$ corresponding to a minimally degenerate system of equations as follows. Taking a generic tuple $(f_2,\ldots,f_n)\in\C^{A'}$, the curve $M=m_{A_1}\{f_2=\ldots=f_n=0\}$ is reduced, irreducible and not a line. Thus, a generic tangent hyperplane $\sum_{a\in A_1} c_az_a=0$ to $M$ has a simple tangency and is transversal to $M$ at the other intersection points. Then the system of equations $\sum_{a\in A_1} c_ax^a=f_2(x)=\ldots=f_n(x)=0$ has one root of multiplicity 2, and the other roots are of multiplicity 1.
\end{rem}

\vspace{1ex}

{\bf The proof of Theorem  \ref{thtransp}.} We first need an explicit construction of the exceptional set $B_A$ in the Kouchnirenko--Bernstein theorem \ref{thbernst}.

The restriction of a linear function $v:\R^n\to\R$ to a finite set $A\subset\Z^n$ takes its maximal value at certain points of $A$. The set of all such points will be denoted by $A^v$. For a tuple $A=(A_1,\ldots,A_k)$, denote the tuple $(A_1^v,\ldots,A_k^v)$ by $A^v$, and the naive discriminant of $A^v$ (see Definition \ref{defdiscr3}) by $D_v$. We shall consider $D_v$ as a subset of $\C^A$ in the sense of Remark \ref{remforget}. The set
$$B=\bigcup_{v\in\R^n} D_v\subset\C^A$$
is algebraic, because there are only finitely many distinct algebraic sets among $D_v,v\in\Z^n$. More specifically, write $u\sim v$ if $A^u=A^v$, then this equivalence relation splits $\R^n$ into finitely many relatively open polyhedral cones. These cones form a fan $\Sigma$ (see e.g. \cite{fulton}), and $A^v$ and $D_v$ depend only on the cone $C\in\Sigma$ containing $v$. So we shall also denote $A^v$ and $D_v$ by $A^C$ and $D_C$ respectively.

We claim that the set $B$ can be taken as the exceptional set $B_A$ in Theorem \ref{thbernst} in the following strong sense. Denote the incidence set $$\{(x,f)\,|\, f(x)=0\}\subset\CC^n\times\C^A$$ by $E$ and its projection to $\C^A$ by $\pi$.

\begin{theor}[(refined Theorem \ref{thtransp})] \label{luseful} Let the tuple $A=(A_1,\ldots,A_n)$, the set $B$ and the projection $\pi$ be as above with $k=n$.

1) The projection $\pi$ is a covering outside the set $B$. In particular, every $f\in\C^A\setminus B$ has exactly $\MV(A)$ roots, and the group $G_A$ is the monodromy group of this covering.

2) If $A$ is reduced and dual effective, then, for a generic $f\in D_A$: 

-- the system $f=0$ has a unique singular root $x\in\CC^n$, and its multiplicity is 2;

-- we have $f\notin D_v$ for every non-zero $v:\R^n\to\R$.

3) For such $f$, let $F:(\C,0)\to (\C^A,f)$ be a germ of a smooth curve transversal to $D_A$. Then the monodromy of the covering from Part (1) along the loop $F\bigl(\varepsilon\exp(2\pi i t)\bigr)$ for small $\varepsilon>0$ is a transposition.
\end{theor}
\begin{rem} \label{forexa1} Instead of assuming dual effectiveness in Part 2, it is enough to assume that the Cayley configuration $A_{\{1,\ldots,n\}}$ is dual effective, and then a small loop around the Cayley discriminant $D_{A_{\{1,\ldots,n\}}}$ still gives a transposition, see Theorem \ref{luseful2} below for this and some other generalizations.
\end{rem}
{\sc Proof of Part 2.} The first statement follows from Corollary \ref{lirrdef}, the second one from Corollary \ref{cordiff} applied to the Cayley discriminant $D_{A_{\{1,\ldots,n\}}}$.

\vspace{2ex}

{\sc Proof of Parts 1 and 3.}  Choose a unimodular simplicial fan $\Sigma'$, subdividing $\Sigma$ (see \cite{kkms} for its existence), and consider the corresponding smooth toric variety $X\supset\CC^n$. Every cone $C\in\Sigma'$ corresponds to an orbit $O_C\subset X$, and, for $C\ne\{0\}$, the closure of the incidence set $E$ in $X\times C^A$ contains a point of the form $(x,f)\in O_C\times\C^A$ only if $f\in D_C$. In particular, if $f\notin D_v$ for every non-zero $v:\R^n\to\R$, then, for a small neighborhood $U\ni f$, its preimage $V=\pi^{-1}(U)$ is disjoint from the orbits $O_C,\, C\ne\{0\}$, i.e. the restriction $\pi:V\to U$ is proper. Now consider two cases, corresponding to the setting of Part 1 and Part 3 respectively: $f\notin D_A$ and $f\in D_A$.

If $f\notin D_A$, then the restriction $\pi:V\to U$ also has no critical points (this claim makes sense, because $E$ is smooth), so it is a trivial covering, and Part 1 is proved.

If $f\in D_A$ has a unique multiple root $x$, and this root has multiplicity 2, then the local degrees of $\pi$ at the point $(x,f)$ and at the other points of the fiber $\pi^{-1}(f)$ equal two and one respectively. Thus $\pi$ has an $\mathcal{A}_1$ singularity at $(x,f)$ and no singularities at other points of the fiber $\pi^{-1}(f)$, i.e. $\pi(z_1,z_2,\ldots,z_N)=(z_1^2,z_2,\ldots,z_N)$ in suitable local coordinates $(z_1,\ldots,z_N)$ on $T$ near $(x,f)$. In particular, the monodromy along a small loop around the origin in the complex line $z_2=\ldots=z_N=0$ is a transposition. $\quad\square$

\vspace{1ex}

{\bf Monodromy of non-square systems of equations.} We outline a generalization of Theorem \ref{luseful} to some reducible tuples $A$ and to the case $k<n$ in order to clarify what happens in examples similar to \ref{excayley} and what could be a natural counterpart of the topic of this paper for non-square systems of equations.

\begin{theor}\label{luseful2} 
Let $A=(A_1,\ldots,A_k)$, $B$ and $\pi$ be as above with arbitrary $k\leqslant n$.

1) The projection $\pi$ in a locally trivial fibration outside the set $B$. Moreover, $B$ is the minimal closed set with this property. In particular, every loop in the complement to $B$ gives rise to the monodromy automorphism in the cohomology $H$ of the fiber of this fibration. 

2)The set $B$ is a hypersurface unless $m+2$ of the sets in the tuple $A$ can be shifted to the same $m$-dimensional plane (in which case the aforementioned fibration is empty).

3) For a reduced tuple $A$, whose Cayley discriminant (Definition \ref{defdiscr3}) is a hypersurface, and for a generic $f$ in the Cayley discriminant: 

-- the system $f=0$ has a unique singular root $x\in\CC^n$, and this singular root is non-degenerate;

-- we have $f\notin D_v$ for every non-zero $v:\R^n\to\R$.

4) For such $f$, let $F:(\C,0)\to (\C^A,f)$ be a germ of a smooth curve transversal to $D_A$. Then the $\zeta$-function of the monodromy transformation from Part 1, corresponding to the loop $F\bigl(\varepsilon\exp(2\pi i t)\bigr)$ for small $\varepsilon>0$ has the form $t^2-1$.
\end{theor}
Parts 1 and 2 follow from Theorems 1.1 and 1.4 in \cite{eadv}. The first statement of Part 3 follows from the fact that the Cayley discriminant is a component of multiplicity 1 in the Euler discriminant $E_A$, see Proposition 1.11 in \cite{eadv}. (This in particular works for $k=n$, but we preferred to give a more straightforward proof of Theorem \ref{luseful} in that case.) The rest is proved in the same way as for $k=n$ in Lemma \ref{luseful}.

\begin{rem} In particular, the correspondence from Theorem \ref{luseful2}.1 maps the fundamental group of the complement of $B$ to the group $GL(H)$. The image $G_A$ is the monodromy group of the (non-square) system of equations supported at the tuple $A$. The results of the present paper give some hope that $G_A$ can be quite explicitly described in terms of $A$ at least for reduced irreducible $A$. This important study has been recently initiated in the simplest non-square case, i.e. $(k,n)=(1,2)$, see \cite{lang1}, \cite{lang2}, \cite{salter}.
\end{rem}

\section{Double transitivity of monodromy} \label{S2tr}

Consider a morphism $\pi$ of an algebraic set $E$ to an irreducible algebraic set $C$ as an {\it abstract enumerative problem}: regard a point $z\in C$ as an incidence condition, and the points of its fiber $\pi^{-1}(z)$ as the solutions of the enumerative problem with a given incidence condition. The enumerative problem is said to be {\it well posed} if its generic fiber is finite. In this case, there exists a Zariski open set $U\subset C$ such that $\pi$ is a covering over $U$. The monodromy group of this covering does not depend on the choice of $U$ and is called the monodromy group of the enumerative problem.
\begin{exa}\label{ourcase} The enumerative problem of the present paper falls into this scheme, if we define $$E=\{(x,f)\,|\, f(x)=0\}\subset\CC^n\times\C^A$$ and denote the projection of $E$ to $C=\C^A$ by $\pi$. For every tuple $A=(A_1,\ldots,A_n)$, it is well posed by Theorem \ref{luseful}.1.
\end{exa}

Let us recall a classically known geometric criterion for the double transitivity of the monodromy of the abstract enumerative problem $\pi:E\to C$. Although its versions can be found in \cite{schsot} and other relevant works, we shall recall the proof to keep the story self-contained. Consider the fiber square $$E_2=\{(x,y)\,|\, \pi(x)=\pi(y)\}\subset E^2$$ and its projection $\pi_2:E_2\to C$, sending $(x,y)$ to $\pi(x)=\pi(y)\in C$.
If the enumerative problem is well posed, i.e., for a certain Zariski open $U\subset C$, its preimage $V=\pi^{-1}(U)$ defines a covering $\pi:V\to U$, then the fiber square $V_2=\pi_2^{-1}(U)$ also defines a covering $\pi_2:V_2\to U$. Note that the diagonal $D=\{(x,x)\,|\, x\in U\}\subset V_2$ is an irreducible component of $V_2$.
\begin{theor} The monodromy of the well posed enumerative problem $\pi:E\to C$ is doubly transitive if and only if $V_2$ has at most one irreducible component different from $D$.
\end{theor}
{\sc Proof.} If $D=V_2$, then the monodromy is trivial. Otherwise, let $F$ be the second component of $V_2$. In order to prove the double transitivity, we should take two pairs of distinct points $(x,y)$ and $(x',y')$ in the fiber $\pi^{-1}(z)$ of a point $z\in U$ and construct a loop in $U$ such that the monodromy along this loop sends $x$ to $x'$ and $y$ to $y'$. Since neither $(x,y)$ nor $(x',y')$ is contained in $D$, both of them are contained in $F$. Since $F$ is irreducible, these two points can be connected with a path $\gamma$. Then $\pi_2(\gamma)$ is the sought loop. $\quad\square$
\begin{sledst}\label{sqcrit} Let $\pi:E\to C$ be a well posed enumerative problem. If at most one irreducible component  of $E_2$ besides the diagonal $D$ has the same dimension as $D$, then the monodromy is doubly transitive.
\end{sledst}

\vspace{1ex}

{\sc Proof of Theorem \ref{th2tr}.}  
The idea is to apply Corollary \ref{sqcrit} to the setting of Example \ref{ourcase}. In this case we have
$$S=\CC^n\times\CC^n,$$ $$E_2=\{(x,y,f)\,|\,f(x)=f(y)=0\}\subset S\times\C^A.$$ 
In order to prove that $G_A$ is doubly transitive, it is enough to prove that $E_2$ has at most one more irreducible component $F$ of dimension $N=\dim \C^A$. We shall prove it by counting the dimension of fibers of the projection $p:E_2\to S$. Every such fiber is a vector subspace of $\C^A$, but different fibers may have different dimension. Namely, assuming for convenience without loss of generality that every $A_i$ contains 0, the fiber $p^{-1}(x,y)$ is given in $\C^A$ by  

$$2n-d_{x,y}\eqno{(*)}$$
independent linear equations, where $d_{x,y}$ is the number of $A_i$'s such that $x^a=y^a$ for all $a\in A_i$. Indeed, since $0\in A_i$, the linear equations $f_i(x)=f_i(y)=0$ on the element $f=(f_1,\ldots,f_i,\ldots,f_n)\in\C^A$ are dependent if and only if they coincide and if and only if $x^a=y^a$ for all $a\in A_i$, so $(*)$ follows. 

This implies that $\dim p^{-1}(x,y)$ is the same for all $(x,y)$ in the set $U_L$, defined as follows:

$$V_L=\{(x,y)\,|\, x^a=y^a \mbox{ for all } a\in L\}\subset S \mbox{ for a sublattice } L\subset\Z^n,$$ 
$$L_I\subset\Z^n \mbox{ is the sublattice generated by } A_i,\, i\in I,$$
$$U_L=V_L\setminus\bigcup_{L_I\supsetneq L} V_{L_I}.$$ 

Namely, if $(x,y)\in U_L$, then, by $(*)$, the fiber $p^{-1}(x,y)$ is given in $\C^A$ by $$2n-d_L\eqno{(**)}$$ independent linear equations, where $d_L$ is the number of $A_i$'s contained in $L$. 

Therefore, denoting the preimage of $U_L$ in $E_2$ by $E_L$, we conclude by $(**)$ that $p:E_L\to U_L$ is a vector bundle of rank $N-2n+d_L,\, N=\dim\C^A$. Moreover, since $\dim U_L=2n-\dim L$, we conclude that $\dim E_L=N+d_L-\dim L$. 

Since the tuple $A$ is reduced and irreducible, we have $$\dim E_L=N+d_L-\dim L<N$$ unless $L=\Z^n$, or $L$ contains no $A_i$'s at all. In the latter cases, $E_L$ equals the diagonal $D\subset E_2$ or one more $N$-dimensional subset $F\subset E_2$ (independent of $L$) respectively. Since $E_2$ is covered by $E_L$'s as $L$ runs over all sublattices, we have proved that it has two $N$-dimensional components, so that Corollary \ref{sqcrit} applies. $\quad\square$

\begin{acknowledgements} I am grateful to Christopher Borger and Benjamin Nill, whose proof \cite{bn} of a conjecture from \cite{cdd} contributed to working out the present approach to the conjecture from \cite{abelruff} and to Yuri Burman and the referee for valuable remarks.\end{acknowledgements}

\medskip

\begin{center}
\section{Addendum}
{\small\it July 7, 2020}
\end{center}

Apparently, I cannot prove at the request of the readers the following step in the proof of Theorem 1.11.

\begin{conj} Given an irreducible tuple $A$ of finite sets is $\Z^n$, 
then every linearly independent tuple $B$ that entirely consists of sets from $A$, can be obtained from $A$ by taking a sequence of AF-descendants $A',A'',\ldots,A^{(k)}=B$.
\end{conj}


This statement implies an estimate 
$$\MV B\leqslant (\MV A)^{2^N}, \eqno{(**)}$$
which leads to the proof of Theorem 1.11.
\begin{rem} 
Conjecture 5.1 is checked directly for small dimensions $n\le 4$, in particular the inequality $(**)$ remains valid in such dimensions.
\end{rem}

The aim of Addendum is to prove the following slightly different version of the estimate $(**)$, preserving the rest of the proof of Theorem 1.11 and all other parts of the original paper.

For $\alpha\in\Z^n$ with non-negative entries, denote the set of $i$ such that $\alpha_i>0$ by $\supp\alpha$, the sum of $\alpha_i$'s by $|\alpha|$, and the tuple containing $\alpha_i$ copies of the set $A_i$ from the tuple $A$ for every $i$ by $\alpha A$. We aim at proving the following.
\begin{theor}
If $A$ is an irreducible tuple of sets in $\Z^n$ of mixed volume $V>0$, and $|\alpha|=n$, then $\MV\alpha A\leqslant (2V)^{2^{n-|\supp\alpha|}}$.
\end{theor}
\begin{rem}
For a tuple of full-dimensional sets (i.e. the ones not contained in an affine hyperplane) a much better estimate is now established by Averkov, Borger ans Soprounov\footnote[1]{G. Averkov, C. Borger, I. Soprunov, {\it Inequalities between mixed volumes of convex bodies: volume bounds for the Minkowski sum}, 	arXiv:2002.03065}. I have been already aware of their work when writing this addendum, and this definitely helped me to come up with using the so called square inequality\footnote[7]{S. Brazitikos, A. Giannopoulos, and D.-M. Liakopoulos, {\it Uniform cover inequalities for the volume of coordinate sections and projections of convex bodies}, Adv. Geom. 18 (2018), no. 3,
345–354, arXiv:1606.03779; c.f. the proof of Lemma 7.4.1 in R. Schneider, {\it Convex bodies: the Brunn-Minkowski theory, expanded ed.}, Encyclopedia of Mathematics
and its Applications, vol. 151, Cambridge University Press, Cambridge, 2014} here
(although I was aware of this inequality by the time of writing my original  paper):
$$\MV(A_1,A_2,A')\MV(A_3,A_3,A')\leqslant 2\MV(A_1,A_3,A')\MV(A_2,A_3,A'),\; A'=(A_3,\ldots,A_n).\eqno{\rm (SQ)}$$
\end{rem}

The rest of the text is devoted to the proof of Theorem 5.3. 

\begin{defin}
1) The codimension $\codim_\alpha I$ of a set $I\subset\supp\alpha$ is $\sum_{i\in I}\alpha_i-\dim\sum_{i\in I}A_i$.

2) A tuple $\alpha$ with $|\alpha|=n$ is said to be linearly independent, if $\supp\alpha$ has no subsets of positive codimension.
\end{defin}
The following is well known (see e.g. [Kh78]).
\begin{lemma}
$\MV\alpha A>0$ if and only if $\alpha$ is linearly independent.
\end{lemma}
Thus, it is enough to prove the theorem for linearly independent tuples. For them, the theorem will be proved by induction on $n-|\supp\alpha|$ (starting from the obvious base $|\supp\alpha|=n$, i.e. $\alpha A=A$).
The inductive step will reduce the case of an arbitrary linearly independent $\alpha$ to the case of one or two tuples with a larger (by one) support. 

To find such larger-supported tuples, we need the following key lemma. Let $e_1,\ldots,e_n$ be the standard basis of $\Z^n$. 

\begin{lemma}\label{keyl}
Assume that $A$ is an irreducible tuple. Then, for every linearly independent $\alpha$ such that $|\alpha|=n$ and $n\not\in\supp\alpha$, there exist $i$ and $j\in\supp\alpha$ (possibly equal) such that $\alpha_i>1, \alpha_j>1$ and $\alpha-e_i-e_j+2e_n$ is linearly independent.
\end{lemma}
The proof is given after the end of the proof of Theorem 5.3. 

The two preceding lemmas prove the first inequality in the following chain, and the second inequality in this chain is ${\rm (SQ)}$: $$2\MV\alpha A\leqslant 2\MV\alpha A\cdot \MV(\alpha-e_i-e_j+2e_n)A\leqslant (2\MV(\alpha-e_i+e_n)A)\cdot (2\MV(\alpha-e_j+e_n)A).\eqno{(1)}$$

Now the statement of the theorem for a linearly independent tuple $\alpha$ follows from the same for $\alpha-e_i+e_n$ and $\alpha-e_j+e_n$ by the inductive hypothesis and $(1)$. Theorem 5.3 is proved.

It remains to prove Lemma \ref{keyl}. Denote $\alpha+2e_n$ by $\tilde\alpha$ (note that $|\tilde\alpha|=n+2$).

\begin{defin}
An obstacle is a subset of $\supp\tilde\alpha$ of positive codimension.
\end{defin}
Since $\alpha$ is linearly independent, every obstacle $I$ contains $n$, and its codimension $\codim_{\tilde\alpha} I$ equals either 1 or 2.

Since the codimension is supermodular (see [St94]) in the sense that
$$\codim_{\tilde\alpha} I+\codim_{\tilde\alpha} J\le\codim_{\tilde\alpha} I\cap J+\codim_{\tilde\alpha} I\cup J,\eqno{(2)}$$
the linear independence of $\alpha$ implies the following three facts.

I. Among all codimension 2 obstacles, there exists the (unique) $I_0$ contained in all the others. This is because the intersection of codimension 2 obstacles is a codimension 2 obstacle by $(2)$, and $\supp\tilde\alpha$ itself is a codimension 2 obstacle.

II. Every minimal (by inclusion) obstacle is contained in $I_0$. This is because the intersection of codimension 1 and codimension 2 obstacles is an obstacle by $(2)$.

III. The union of every two (distinct) minimal obstacles equals $I_0$. This is because otherwise their intersection would be an obstacle by $(2)$.

Denoting the codimension 1 minimal obstacles by $I_1,\ldots,I_q$, we summarize (I-III) as follows:

$$I_k\subset I_0 \mbox{ and } I_k\cup I_m=I_0 \mbox{ for all } k\ne m.$$

If there exists $i_0\in I_0\cap\ldots\cap I_q$ such that $\alpha_{i_0}>1$, then we prove Leamma \ref{keyl} by setting $(i,j):=(i_0,i_0)$.

Even if such $i_0$ does not exist, the irreducibility of $A$ implies that every obstacle $I_s$ contains its own element $i_s\ne n$ such that $\alpha_{i_s}>1$.

Notice that $i_s$ is not contained in some $I_t$ (otherwise it could be taken as $i_0$ above). Then $i_s$ is contained in all other $I_r,\,r\ne t$, because $I_r\cup I_t=I_0$. Thus we can prove Lemma \ref{keyl} in this remaining case by setting $(i,j):=(i_s,i_t)$.

\end{document}